\documentclass[letterpaper, 10 pt, conference]{ieeeconf}  

\IEEEoverridecommandlockouts                              
\usepackage[utf8]{inputenc}
\usepackage[T1]{fontenc}
\usepackage{amsmath}
\usepackage{mathrsfs}
\usepackage{amsfonts}
\usepackage{algorithm}
\usepackage{algorithmicx}
\usepackage{algpseudocode}

\usepackage{color}
\usepackage{hyperref}

 \usepackage{graphicx}

\title{\LARGE \bf
Lifted contact dynamics for efficient optimal control of \\
rigid body systems with contacts 
}

\author{Sotaro Katayama$^{1}$ and Toshiyuki Ohtsuka$^{1}$
\thanks{$^{1}$S. Katayama and T. Ohtsuka are with the Department of System Science, Graduate School of Informatics, Kyoto University, Kyoto, Japan
        {\tt\small katayama@sys.i.kyoto-u.ac.jp}, 
        {\tt\small ohtsuka@i.kyoto-u.ac.jp}}%
}

\begin{document}

\onecolumn
\noindent
© 2022 IEEE. Personal use of this material is permitted. Permission from IEEE must be obtained for all other uses, in any current or future media, including reprinting/republishing this material for advertising or promotional purposes, creating new collective works, for resale or redistribution to servers or lists, or reuse of any copyrighted component of this work in other works.

\hspace{1cm}

\noindent
\textbf{Published article:} \\ 
\noindent
S. Katayama and T. Ohtsuka, ``Lifted contact dynamics for efficient optimal control of rigid body systems with contacts,'' 2022 IEEE/RSJ International Conference on Intelligent Robots and Systems (IROS2022), 2022, pp. 8879--8886.

\twocolumn

\maketitle
\thispagestyle{empty}
\pagestyle{empty}

\begin{abstract}

We propose a novel and efficient lifting approach for the optimal control of rigid-body systems with contacts to improve the convergence properties of Newton-type methods.
To relax the high nonlinearity, we consider the state, acceleration, contact forces, and control input torques, as optimization variables and the inverse dynamics and acceleration constraints on the contact frames as equality constraints.
We eliminate the update of the acceleration, contact forces, and their dual variables from the linear equation to be solved in each Newton-type iteration in an efficient manner. 
As a result, the computational cost per Newton-type iteration is almost identical to that of the conventional non-lifted Newton-type iteration that embeds contact dynamics in the state equation. 
We conducted numerical experiments on the whole-body optimal control of various quadrupedal gaits subject to the friction cone constraints considered in interior-point methods and demonstrated that the proposed method can significantly increase the convergence speed to more than twice that of the conventional non-lifted approach.

\end{abstract}

\section{Introduction}
Model predictive control (MPC) is a promising framework for reactive motion planning and control of rigid-body systems that make rigid contacts with the environment, e.g., legged robots and robot manipulators. 
A remaining critical constraint of MPC is the computational time required to solve optimal control problems (OCPs) with limited computational resources.
Here, we introduce computational issues for two classes of solution approaches of OCPs of robotic systems: contact-implicit and event-driven approaches.

Contact-implicit approaches aim to determine the optimal trajectory without specifying the contact sequence in advance. 
The simplest approach is to approximate the contact forces using spring-damper systems \cite{bib:legged2};  
however, it lacks accuracy and involves stiff optimization problems.
The same problem typically occurs in other smooth soft-contact models such as in \cite{bib:todorovmodel:apply}.
More accurate approaches are based on the mathematical programs with complementarity constraints (MPCCs) \cite{bib:implicitContact} or bilevel optimization (BO) problems \cite{bib:implicitHardContact}.
However, the MPCCs inherently lack linear independence constraint qualification, which causes practical and theoretical convergence issues \cite{bib:MPCCLimits}.
Moreover, the convergence analysis of BO is yet to be discussed \cite{bib:implicitHardContact}, which renders the practical applications of BO difficult.

In contrast, event-driven approaches formulate the OCPs under predefined contact sequences and model the impacts between the system and the environment explicitly via Newton's law of impact.
Therefore, event-driven approaches are more accurate and the continuous-time OCP can be discretized with coarser grids than in contact-implicit methods, which means that they can be faster.  
Moreover, they consist of only smooth optimization problems; therefore, we can solve them efficiently using the off-the-shelf Newton-type methods for smooth OCPs without being overly cautious with the LICQ. 
Although event-driven approaches cannot optimize contact sequences as contact-implicit methods can, they can optimize the instants of the impacts under a given contact sequence, as studied in \cite{bib:DDP:jumpRobot}. 

\begin{figure}
\centering
\includegraphics[scale=0.335]{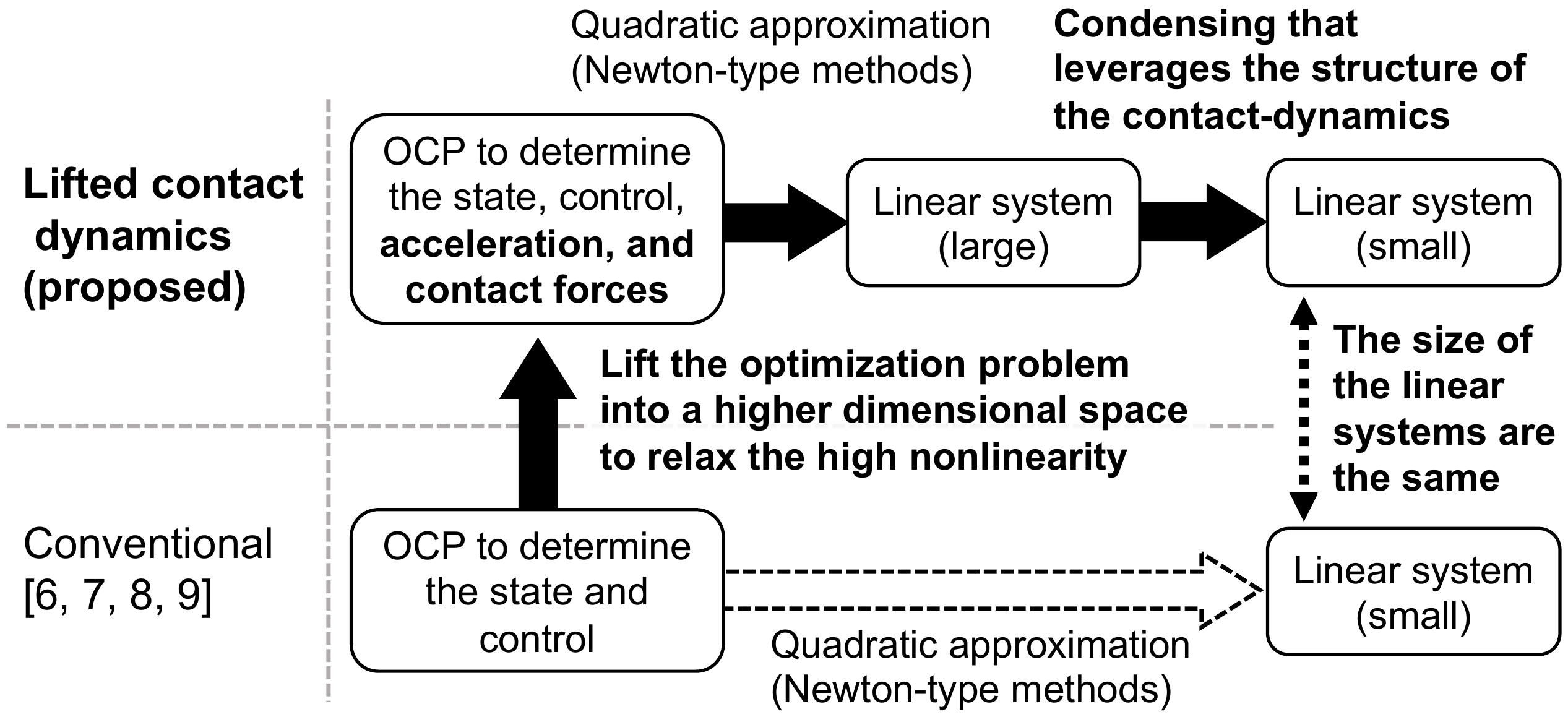}
\caption{
Conceptual diagram of the proposed lifted contact dynamics scheme for efficient optimal control of rigid body systems with contacts.
The conventional optimal control problem (OCP) is lifted into a higher dimensional one to relax the high nonlinearity. 
The Newton-type method is carried out efficiently by leveraging the structure of the contact dynamics.
}
\label{fig:conceptual}
\end{figure}

The contact-consistent forward dynamics, a calculation of the acceleration and contact forces using the given configuration, velocity, and torques, which we call ``contact dynamics,” has been successfully used with Newton-type methods to efficiently solve event-driven OCPs \cite{bib:DDP:jumpRobot, bib:humanLikeRunnig, bib:DDPContact, bib:crocoddyl}.
The contact dynamics formulation is utilized in \cite{bib:humanLikeRunnig} with the direct multiple shooting method (DMS) \cite{bib:condensing}.
\cite{bib:DDPContact} and \cite{bib:DDP:jumpRobot} solve the OCPs using differential dynamic programming (DDP) \cite{bib:DDPBook}.
\cite{bib:crocoddyl} improved the approach of \cite{bib:DDPContact} using the efficient analytical derivatives proposed in \cite{bib:analytical:rbd} to compute the sensitivities of contact dynamics.
However, such an approach still contains high nonlinearity, which results in slow convergence, particularly when there are costs and constraints on the contact forces, e.g., friction cone constraints.

The lifted Newton method is a promising option to relax such high nonlinearity. It involves adding intermediate variables and lifting the optimization problem into a higher-dimensional one \cite{bib:liftedNewton}.
For example, it is well known in the context of numerical optimal control and MPC that the multiple shooting method, which considers the state and control input as the optimization variables, has preferred convergence properties over the single shooting method, which considers only the control input as the optimization variable \cite{bib:condensing}.
However, to the best of our knowledge, no lifting method has previously been used to efficiently treat the contact forces in event-driven OCPs.
For example, in our previous study \cite{bib:InvDynOCP}, we utilized a certain lifted formulation; however, it may be inefficient when the number of contacts is significant. It is efficient only for systems without contacts or with a few contacts by leveraging inverse dynamics computation.
Moreover, since the focus of the previous study was on the computational time of rigid body systems without contacts, no discussion or numerical investigation regarding convergence properties under contacts were conducted.

In this paper, we propose \textit{lifted contact dynamics}, a novel and efficient lifting approach for the optimal control of rigid-body systems with contacts to improve the convergence properties of Newton-type methods. 
Fig. \ref{fig:conceptual} illustrates the proposed framework.
To relax the high nonlinearity, we consider the state, acceleration, contact forces, and control input torques, as the optimization variables and the inverse dynamics and acceleration constraints on the contact frames as equality constraints.
We eliminate the updating of the acceleration, contact forces, and their dual variables from the linear equation to be solved in each Newton-type iteration in an efficient manner.
As a result, the computational cost per Newton-type iteration is almost identical to that of the conventional non-lifted Newton-type iteration that embeds contact dynamics in the state equation.
This aspect distinguishes this study from our previous inverse dynamics-based algorithm \cite{bib:InvDynOCP}, which was inefficient for high-dimensional contact constraints.
We also propose a lifting method for impulse dynamics, which was not considered in our previous study \cite{bib:InvDynOCP}.
We conducted numerical experiments on the whole-body optimal control of various quadrupedal gaits subject to the friction cone constraints considered in the interior-point methods, while \cite{bib:InvDynOCP} treated much simpler examples.
The experiments demonstrated that the proposed method can significantly increase the convergence speed to more than twice that of conventional approaches based on the non-lifted contact dynamics.
The contributions of this paper are then summarized as follows:
\begin{itemize}
    \item A lifted formulation for OCP of robotic systems with contacts to relax the high nonlinearity.
    \item Efficient ``condensing'' algorithms to enable as fast as Newton step computation as the non-lifted counterpart. 
    \item Numerical studies on the practical quadrupedal locomotion problems.
\end{itemize}

The remainder of this paper is organized as follows. In Section \ref{section:formulation}, we review the contact and impulse dynamics and conventional formulations of event-driven OCPs.
In Section \ref{section:lifted}, we introduce the lifted contact dynamics, a Newton-type method that efficiently condenses the linear equations to be solved in each Newton iteration. 
Section \ref{section:expriment} compares the proposed method with existing approaches based on the non-lifted contact dynamics and demonstrates its effectiveness in terms of the convergence properties. 
Section \ref{section:conclu} concludes the paper and outlines future research directions.

\textit{Notation:} We denote the partial derivatives of a differentiable function with certain variables using a function with subscripts; i.e., $f_x (x)$ denotes $\frac{\partial f}{\partial x} (x)$ and $g_{x y} (x, y)$ denotes $\frac{\partial^2 g}{\partial x \partial y} (x, y)$. 
We denote an $n \times n$ identity matrix as $I_{n}$.

\section{Overview of Optimal Control Problems of Rigid-Body Systems with Contacts}\label{section:formulation}
\subsection{Contact dynamics}
First, we review the components to formulate event-driven direct OCPs, e.g., the contact dynamics and discretized state equations.
Let $Q$ be the configuration manifold of the rigid-body system. Let $q \in Q$, $v \in \mathbb{R}^n$, $a \in \mathbb{R}^n$, $f \in \mathbb{R}^{n_f}$, and $u \in \mathbb{R}^{n_a}$ be the configuration, generalized velocity, acceleration, stack of the contact forces, and torques of the actuated joints, respectively. The equation of motion of the rigid-body system is expressed as
\begin{equation}\label{eq:equationOfMotion}
    M(q) a + h(q, v) - J^{\rm T} (q) f = S^{\rm T} u,
\end{equation}
where $M (q) \in \mathbb{R}^{n \times n}$ denotes the inertia matrix, $h (q, v) \allowbreak \in \mathbb{R}^{n}$ encompasses the Coriolis, centrifugal, and gravitational terms, $J (q) \in \mathbb{R}^{n_f \times n}$ denotes the stack of the contact Jacobians, and $S \in \mathbb{R}^{n_a \times n}$ denotes the selection matrix. 
The evolution of the state $\begin{bmatrix} q ^{\rm T} & v ^{\rm T} \end{bmatrix} ^{\rm T}$ with time step $\Delta \tau > 0$ is expressed as 
\begin{equation}\label{eq:qvEvolution}
\begin{bmatrix}
    q^+ \\
    v^+
\end{bmatrix}
= \begin{bmatrix}
    q \oplus v \Delta \tau \\ 
    v + a \Delta \tau
\end{bmatrix},
\end{equation}
where $\oplus$ denotes the increment operator on the configuration manifold $Q$ \cite{bib:microLieTheory}.
We formulate the OCPs based on the DMS and then consider an equivalent equality constraint:
\begin{equation}\label{eq:qvDiff}
\begin{bmatrix}
    \delta (q, q ^+) + v \Delta \tau \\
    v - v ^+ + a \Delta \tau
\end{bmatrix} = 0,
\end{equation}
where $\delta (q_1, q_2) := q_1 \ominus q_2 \in \mathbb{R}^n$, and $\ominus$ denotes the subtraction operator between the two configurations \cite{bib:microLieTheory}.
The system also satisfies the contact constraints of the form 
\begin{equation}\label{eq:p}
    {\bf p} (q) = 0,
\end{equation}
where ${\bf p} (q) \in \mathbb{R}^{n_f}$ is the stack of the positions of the contact frames.
In event-driven OCPs, instead of considering (\ref{eq:p}) over a time interval, we typically constraint on the acceleration of the contact frames over the time interval that is generalized into the form of Baumgarte's stabilization method \cite{bib:baumgarte}:
\begin{equation}\label{eq:a}
    {\bf a} (q, v, a) := \ddot{\bf p} + 2 \alpha \dot{\bf p} + \beta ^2 {\bf p} = J (q) a + {\bf b} (q, v),
\end{equation}
where $\alpha$ and $\beta$ are weight parameters, and we define
\begin{equation}\label{eq:b}
  {\bf b} (q, v) := \dot{J} (q, v) v + 2 \alpha J (q) v  + \beta^2 {\bf p} (q).
\end{equation}
By setting $\alpha = \beta > 0$, we can stabilize the violation of the original constraint (\ref{eq:p}) over the time interval provided that (\ref{eq:p}) and the equality constraint on the contact velocity, 
\begin{equation}\label{eq:v}
    \dot{\bf p} (q, v) = J(q) v = 0,
\end{equation}
are satisfied at a point on the interval \cite{bib:baumgarteParameters}.
Note that (\ref{eq:a}) is reduced to the equality constraint on the acceleration of the contact frames with $\alpha = \beta = 0$.
By combining (\ref{eq:equationOfMotion}) and (\ref{eq:a}), we obtain the contact dynamics, i.e., the contact-consistent forward dynamics: 
\begin{equation}\label{eq:contactDynamics}
    \begin{bmatrix}
        M (q) & J^{\rm T} (q) \\ 
        J (q) & O
    \end{bmatrix}
    \begin{bmatrix}
        a \\ 
        - f
    \end{bmatrix}
    = 
    \begin{bmatrix}
        S ^{\rm T} u - h (q, v) \\
        - {\bf b} (q, v)
    \end{bmatrix}.
\end{equation}
In the conventional approaches \cite{bib:DDP:jumpRobot, bib:humanLikeRunnig, bib:DDPContact, bib:crocoddyl}, we eliminate $a$ and $f$ as functions of $q$, $v$, and $u$ from the OCP using (\ref{eq:contactDynamics}).
For example, we substitute $a$ with (\ref{eq:contactDynamics}) in (\ref{eq:qvEvolution}) or (\ref{eq:qvDiff}) and consider these equations to be the discretized state equation.
We can also eliminate $a$ and $f$ from the cost function and the constraints, e.g., the friction cone constraints, using (\ref{eq:contactDynamics}).

\subsection{Impulse dynamics}
Similar to the contact dynamics, the equation of Newton's law of impact of the systems is expressed as 
\begin{equation}\label{eq:equationOfImpulse}
    M (q) \delta v - J ^{\rm T} (q) \Lambda = 0,
\end{equation}
where $\delta v \in \mathbb{R}^n$ denotes the impulse change in the generalized velocity, and $\Lambda \in \mathbb{R}^{n_f}$ denotes the stack of the impact forces.
The evolution of the state between the impulse and its equivalent equality constraint considered in the DMS are expressed as
\begin{equation}\label{eq:qvEvolutionImpulse}
\begin{bmatrix}
    q^+ \\
    v^+
\end{bmatrix}
= \begin{bmatrix}
    q \\ 
    v + \delta v
\end{bmatrix}
\end{equation}
and
\begin{equation}\label{eq:qvDiffImpulse}
\begin{bmatrix}
    \delta (q, q ^+) \\
    v - v ^+ + \delta v
\end{bmatrix} = 0,
\end{equation}
respectively. 
Herein, we assume a completely inelastic collision, which results in the contact velocity constraints of the form (\ref{eq:v}) immediately after the impulse as 
\begin{equation}\label{eq:vplus}
    {\bf v} (q, v, \delta v) := \dot{\bf p} (q, v + \delta v) = J (q) (v + \delta v) = 0.
\end{equation}
The contact position constraint (\ref{eq:p}) for the frames with impacts is also imposed at the impulse instant.
By combining (\ref{eq:equationOfImpulse}) and (\ref{eq:vplus}), we obtain the impulse dynamics:
\begin{equation}\label{eq:impulseDynamics}
    \begin{bmatrix}
        M (q) & J^{\rm T} (q) \\ 
        J (q) & O
    \end{bmatrix}
    \begin{bmatrix}
        \delta v \\ 
        - \Lambda
    \end{bmatrix}
    = 
    \begin{bmatrix}
        0 \\
        - J (q) v
    \end{bmatrix}.
\end{equation}
In the conventional formulation, as well as the contact dynamics, $\delta v$ and $\Lambda$ are eliminated from the OCP as functions of $q$ and $v$ using (\ref{eq:impulseDynamics}), for example, from (\ref{eq:qvEvolutionImpulse}) and (\ref{eq:qvDiffImpulse}).

\subsection{Conventional formulation of optimal control}
We summarize the conventional formulation of the optimal control of rigid-body systems \cite{bib:DDP:jumpRobot, bib:humanLikeRunnig, bib:DDPContact, bib:crocoddyl}.
We introduce $N$ discretization grids. 
We define $\mathcal{J} \subset \left\{ 0,...,N-1 \right\}$ which is the set of the impulse stages and define $\mathcal{I} := \left\{ 0,...,N-1 \right\} \backslash \mathcal{J}$.
For given user-defined terminal cost $V_f (\cdot)$ and stage costs $l_i (\cdot)$, the conventional OCP is summarized as follows: determine $q_0, ..., q_N \in Q$, $v_0, ..., v_N \in \mathbb{R}^n$, and $\left\{u_i \right\}_{i \in \mathcal{I}} \in \mathbb{R}^m$ that minimize a given cost function 
\begin{equation}\label{eq:cost}
    V_f (x_N) + \sum_{i \in \mathcal{I}} l_i (x_i, u_i, a_i, f_i) + \sum_{j \in \mathcal{J}} l_j (x_j, \delta v_i, \Lambda_i),
\end{equation}
where $x_i := \begin{bmatrix} q_i ^{\rm T} & v_i ^{\rm T} \end{bmatrix}^{\rm T}$ is the state, 
subject to the state equation obtained by eliminating $a_i$ from (\ref{eq:qvEvolution}) using (\ref{eq:contactDynamics}), the state equation at the impulse obtained by eliminating $\delta v_j$ from (\ref{eq:qvEvolutionImpulse}) using (\ref{eq:impulseDynamics}), the contact-position constraint (\ref{eq:p}) at the impulse instant, and the other user-defined equality and inequality constraints.
In addition, $a_i$, $f_i$, $\delta v_j$, and $\Lambda_j$ are also eliminated as functions of $q_i$, $v_i$, and $u_i$ from the cost function (\ref{eq:cost}) and the constraints other than the state equations, which increases the nonlinearity.
In \cite{bib:humanLikeRunnig}, this problem was solved using the DMS, that is, $x_0, ..., x_N$ and $\left\{ u_i \right\}_{i \in \mathcal{I}}$ were considered as the optimization variables, the state equations were considered as the equality constraints, and (\ref{eq:cost}) was represented by $x_0, ..., x_N$ and $\left\{ u_i \right\}_{i \in \mathcal{I}}$.
In \cite{bib:DDP:jumpRobot, bib:DDPContact, bib:crocoddyl}, this problem was solved using kinds of the single shooting method, that is, $x_0, ..., x_N$ were further eliminated from the optimization problem using the state equations, and only $\left\{ u_i \right\}_{i \in \mathcal{I}}$ were considered as the decision variables, e.g., (\ref{eq:cost}) was represented only by $\left\{ u_i \right\}_{i \in \mathcal{I}}$.

\section{Lifted Contact Dynamics in Optimal Control}\label{section:lifted}
\subsection{Lifted contact dynamics}
We herein present the lifted contact dynamics to relax the high nonlinearity in OCPs, of which conceptual diagram is shown in Fig. \ref{fig:conceptual}.
Let $i \in \mathcal{I}$, $y_i := ( q_i, v_i, a_i, f_i )$, and $z_i := ( q_i, v_i, a_i )$.
We first augment the control input to convert the system into a fully actuated system in the numerical optimization.
Without loss of generality, 
we assume that $S$ is given by $\begin{bmatrix}
    O & I
\end{bmatrix}^{\rm T}$,
and we define
$\tilde{u}_i := 
\begin{bmatrix}
    u_i ^0 \\
    u_i
\end{bmatrix}$,
where $u_i ^0 \in \mathbb{R}^{n - n_a}$ denotes the virtual torques on the passive joints.
Thus, we can express (\ref{eq:equationOfMotion}) as an equality constraint of the inverse dynamics as follows:
\begin{equation}\label{eq:inverseDynamics}
    {\rm ID} (y_i) - \tilde{u}_i = 0,
\end{equation}
where ${\rm ID} (y_i)$ is defined by the left-hand side of (\ref{eq:equationOfMotion}) and can be computed efficiently using the recursive Newton--Euler algorithm (RNEA) \cite{bib:featherstone}, a fast algorithm for inverse dynamics, with an additional equality constraint
\begin{equation}\label{eq:uPassive}
    u_i ^0 = \bar{S} ^{\rm T} \tilde{u}_i = 0,
    \;\;\; \bar{S} := \begin{bmatrix}
    I & O 
    \end{bmatrix}^{\rm T} \in \mathbb{R}^{n \times n-n_a}.
\end{equation}
We consider $a_i$ and $f_i$ as the optimization variables and (\ref{eq:equationOfMotion}), (\ref{eq:a}), and (\ref{eq:uPassive}) as the equality constraints.
The first-order derivatives of the Lagrangian $\mathcal{L}$ (defined by augmenting the constraints to the cost function (\ref{eq:cost})) with respect to the variables at the stage $i$, that is, the part of the Karush--Kuhn--Tucker (KKT) conditions associated with the stage $i$, are given by (\ref{eq:qvDiff}), (\ref{eq:inverseDynamics}), (\ref{eq:a}), and (\ref{eq:uPassive}), 
\begin{align}\label{eq:lqvi}
    & \begin{bmatrix}
        \mathcal{L}_{q_i} ^{\rm T} \\
        \mathcal{L}_{v_i} ^{\rm T} 
    \end{bmatrix}
    = \begin{bmatrix}
        l_{q_i} ^{\rm T} \\
        l_{v_i} ^{\rm T} 
    \end{bmatrix}
    + \begin{bmatrix}
        \delta_{q_{i}} ^{\rm T} (q_{i}, q_{i+1}) & O \\
        I_n \Delta \tau & I_n 
    \end{bmatrix}
    \begin{bmatrix}
        \lambda_{i+1} \\
        \gamma_{i+1}
    \end{bmatrix} \notag \\
    & + \begin{bmatrix}
        {\rm ID}_{q_i} ^{\rm T} (y_i) & {\bf a}_{q_i} ^{\rm T} (z_i) \\
        {\rm ID}_{v_i} ^{\rm T} (y_i) & {\bf a}_{v_i} ^{\rm T} (z_i)
    \end{bmatrix}
    \begin{bmatrix}
        \beta_i \\
        \mu_i 
    \end{bmatrix} \Delta \tau
    + \begin{bmatrix}
        \delta_{q_{i}} ^{\rm T} (q_{i-1}, q_i) \lambda_{i} \\
        - \gamma_{i}
    \end{bmatrix} = 0,
\end{align}
\begin{align}\label{eq:lafi}
    \begin{bmatrix}
        \mathcal{L}_{a_i} ^{\rm T} \\
        - \mathcal{L}_{f_i} ^{\rm T} 
    \end{bmatrix}
    = & \begin{bmatrix}
        l_{a_i} ^{\rm T} \\
        - l_{f_i} ^{\rm T} 
    \end{bmatrix}
    + \begin{bmatrix}
        \gamma_{i+1} \Delta \tau \\
        0 
    \end{bmatrix} 
    \notag \\ 
    & + \begin{bmatrix}
        {\rm ID}_{a_i} ^{\rm T} (y_i) & {\bf a}_{a_i} ^{\rm T} (z_i) \\
        - {\rm ID}_{f_i} ^{\rm T} (y_i) & - {\bf a}_{f_i} ^{\rm T} (z_i)
    \end{bmatrix}
    \begin{bmatrix}
        \beta_i \\
        \mu_i
    \end{bmatrix} \Delta \tau = 0,
\end{align}
and
\begin{equation}\label{eq:lui}
    \mathcal{L}_{\tilde{u}_i} ^{\rm T} := l_{\tilde{u}_i} ^{\rm T} - \beta_i \Delta \tau + \bar{S} \nu_i \Delta \tau = 0,
\end{equation}
where $\lambda_{i+1}, \gamma_{i+1}, \allowbreak \beta_i \in \mathbb{R}^n$, $\mu_i \in \mathbb{R}^{n_f}$, and $\nu \in \mathbb{R}^{n-m}$ are the Lagrange multipliers with respect to (\ref{eq:qvDiff}), (\ref{eq:inverseDynamics}), (\ref{eq:a}), and (\ref{eq:uPassive}), respectively. 

In the Newton-type methods, these KKT conditions are linearized into linear equations with respect to the Newton steps $\Delta q_i$, $\Delta v_i$, $\Delta a_i$, $\Delta f_i$, $\Delta u_i$, $\Delta u_i ^0$, $\Delta \lambda_{i+1}$, $\Delta \gamma_{i+1}$, $\Delta \beta_i$, $\Delta \mu_i$, and $\Delta \nu_i$.
However, this problem is significantly larger than the conventional OCPs based on the non-lifted contact dynamics that considers only $\Delta q_i$, $\Delta v_i$, $\Delta u_i$, $\Delta \lambda_{i+1}$, and $\Delta \gamma_{i+1}$.
Thus, we propose an efficient condensing method to reduce the size of the linear equation.
We first observe that the equality constraints (\ref{eq:inverseDynamics}) and (\ref{eq:a}) are linearized as 
\begin{align}\label{eq:contactDynamicsLinearized}
    & \begin{bmatrix}
        M (q_i) & J ^{\rm T} (q_i) \\
        J (q_i) & O 
    \end{bmatrix}
    \begin{bmatrix}
        \Delta a_i \\
        - \Delta f_i 
    \end{bmatrix} \notag \\
    & = - \begin{bmatrix}
        {\rm ID}_{q_i} (y_i) & {\rm ID}_{v_i} (y_i) \\
        {\bf a}_{q_i} (z_i) & {\bf a}_{v_i} (z_i)
    \end{bmatrix}
    \begin{bmatrix}
        \Delta q_i \\
        \Delta v_i 
    \end{bmatrix}
    + \begin{bmatrix}
        \Delta \tilde{u}_i \\
        0
    \end{bmatrix}
    - \begin{bmatrix}
        {\rm ID} (y_i) \\
        {\bf a} (z_i)
    \end{bmatrix}.
\end{align}
Furthermore, we always set $u_i ^0 = 0$ and obtain from (\ref{eq:uPassive})
\begin{equation}\label{eq:uPassiveLinearized}
    \Delta u_i ^0 = - u_i ^0 = 0.
\end{equation}
Therefore, we can express $\Delta a_i$ and $\Delta f_i$ using the linear combinations of $\Delta q_i$, $\Delta v_i$, and $\Delta u_i$ using (\ref{eq:contactDynamicsLinearized}) and (\ref{eq:uPassiveLinearized}) if we compute 
\begin{equation}\label{equ:MJtJinv}
\begin{bmatrix}
    M (q_i) & J ^{\rm T} (q_i) \\
    J (q_i) & O 
\end{bmatrix}^{-1}.
\end{equation}
Next, we observe that (\ref{eq:lafi}) and (\ref{eq:lui}) are linearized into 
\begin{align}\label{eq:lafiLin}
    & \begin{bmatrix}
        \mathcal{L}_{a_i} ^{\rm T} \\
        - \mathcal{L}_{f_i} ^{\rm T} 
    \end{bmatrix}
    + \begin{bmatrix}
        \mathcal{L}_{a_i y_i} \\
        - \mathcal{L}_{f_i y_i}
    \end{bmatrix}
    \Delta y_i
    + \begin{bmatrix}
        \Delta \gamma_{i+1} \Delta \tau \\
        O 
    \end{bmatrix}
    \notag \\
    & 
    + \begin{bmatrix}
        M (q_i) & J ^{\rm T} (q_i) \\
        J (q_i) & O
    \end{bmatrix}
    \begin{bmatrix}
        \Delta \beta_i \\
        \Delta \mu_i 
    \end{bmatrix} \Delta \tau = 0 
\end{align}
and 
\begin{equation}\label{eq:luiLin}
    \mathcal{L}_{\tilde{u}_i y_i} ^{\rm T} \Delta y_i - \Delta \beta_i \Delta \tau + \bar{S} \Delta \nu_i \Delta \tau = 0.
\end{equation}
Therefore, we can also express the Newton steps of the dual variables $\Delta \beta_i$, $\Delta \mu_i$, and $\Delta \nu_i$ with the linear combinations of $\Delta q_i$, $\Delta v_i$, $\Delta u_i$, and $\Delta \gamma_{i+1}$ using (\ref{eq:lafiLin}) and (\ref{eq:luiLin}) if we compute (\ref{equ:MJtJinv}).
Subsequently, the linear equations to be solved in the Newton-type iterations at stage $i$ are reduced to that with respect to $\Delta q_i$, $\Delta v_i$, $\Delta u_i$, $\Delta \lambda_{i+1}$, and $\Delta \gamma_{i+1}$, defined as
\begin{align}
    & \begin{bmatrix}
        \tilde{F}_{q, i} \\ 
        \tilde{F}_{v, i}
    \end{bmatrix} 
    + \begin{bmatrix}
        \delta_{q_{i}} (q_{i}, q_{i+1}) & I \Delta \tau \\ 
        \tilde{F}_{vq, i} & \tilde{F}_{vv, i}
    \end{bmatrix} 
    \begin{bmatrix}
        \Delta q_{i} \\ 
        \Delta v_{i}
    \end{bmatrix} 
    + \begin{bmatrix}
        O \\ 
        \tilde{F}_{vu, i}
    \end{bmatrix}
    \Delta u_i
    \notag \\
    &
    + \begin{bmatrix}
        \delta_{q_{i+1}} (q_{i}, q_{i+1}) & O \\ 
        O & - I
    \end{bmatrix}
    \begin{bmatrix}
        \Delta q_{i+1} \\ 
        \Delta v_{i+1}
    \end{bmatrix} = 0
\end{align}
\begin{align}\label{eq:lqviCon}
    & \begin{bmatrix}
        \tilde{\mathcal{L}}_{q_i} ^{\rm T} \\
        \tilde{\mathcal{L}}_{v_i} ^{\rm T} 
    \end{bmatrix}
    + \begin{bmatrix}
        \tilde{\mathcal{L}}_{q_i q_i} &
        \tilde{\mathcal{L}}_{q_i v_i} \\
        \tilde{\mathcal{L}}_{v_i q_i} &
        \tilde{\mathcal{L}}_{v_i v_i} 
    \end{bmatrix}
    \begin{bmatrix}
        \Delta q_i \\
        \Delta v_i
    \end{bmatrix} 
    + \begin{bmatrix}
        \delta_{q_{i}} ^{\rm T} (q_{i-1}, q_i) \Delta \lambda_{i} \\
        - \Delta \gamma_{i}
    \end{bmatrix} 
    \notag \\
    & + \begin{bmatrix}
        \delta_{q_{i}} ^{\rm T} (q_{i}, q_{i+1}) 
        & \tilde{F}_{vq, i} \\
        I_n \Delta \tau 
        & \tilde{F}_{vv, i}
    \end{bmatrix}
    \begin{bmatrix}
        \Delta \lambda_{i+1} \\
        \Delta \gamma_{i+1}
    \end{bmatrix} 
    + \begin{bmatrix}
        \tilde{\mathcal{L}}_{q_i u_i} \\
        \tilde{\mathcal{L}}_{v_i u_i} 
    \end{bmatrix} 
    \Delta u_i = 0,
\end{align}
and
\begin{equation}\label{eq:luiCon}
    \tilde{\mathcal{L}}_{u_i} ^{\rm T} 
    + \begin{bmatrix}
        \tilde{\mathcal{L}}_{q_i u_i} &
        \tilde{\mathcal{L}}_{v_i u_i} 
    \end{bmatrix}
    \begin{bmatrix}
        \Delta q_i \\
        \Delta v_i
    \end{bmatrix} 
    + \tilde{\mathcal{L}}_{u_i u_i} \Delta u_i
    + \tilde{F}_{vu, i} \Delta \gamma_{i+1} = 0,
\end{equation}
where the vectors and matrices with a superscript tilde are derived via simple additions and multiplications of the original Hessians of the Lagrangian, Jacobians of the constraints, and residuals in the KKT conditions.
We omit their definitions here because they are excessively long, but they are not hard to derive.
Note that our previous approach \cite{bib:InvDynOCP} involves solving the linear equation with respect to $\Delta q_i$, $\Delta v_i$, $\Delta a_i$, $\Delta f_i$, $\Delta \mu_i$, and $\Delta \nu_i$; thus, it can be inefficient compared with the proposed method.
After solving the linear equations and obtaining $\Delta q_i$, $\Delta v_i$, $\Delta u_i$, $\Delta \lambda_{i+1}$, and $\Delta \gamma_{i+1}$, we can easily compute $\Delta a_i$, $\Delta f_i$, $\Delta \beta_i$, and $\Delta \mu_i$ using (\ref{eq:contactDynamicsLinearized}) and (\ref{eq:lafiLin}), which is called an expansion procedure.

\subsection{Lifted impulse dynamics}
Next, we present the lifted impulse dynamics.
Let $j \in \mathcal{J}$, $y_j := ( q_j, v_j, \delta v_j, \Lambda_j )$, and $z_j := ( q_j, v_j, \delta v_j )$.
We express (\ref{eq:equationOfImpulse}) as an equality constraint:
\begin{equation}\label{eq:inverseImpulseDynamics}
    \Gamma (y_j) = 0,
\end{equation}
where $r (y_i)$ is defined by the left-hand side of (\ref{eq:equationOfImpulse}).
We consider $\delta v_j$ and $\Lambda_j$ as the optimization variables and (\ref{eq:inverseImpulseDynamics}) and (\ref{eq:v}) as the equality constraints.
The first-order derivatives of the Lagrangian $\mathcal{L}$ with respect to the variables at the stage of impulse $j$, that is, the part of the KKT conditions associated with the impulse stage $j$, are given by (\ref{eq:qvDiffImpulse}), (\ref{eq:inverseImpulseDynamics}), (\ref{eq:v}), 
\begin{align}\label{eq:lqvj}
    & \begin{bmatrix}
        \mathcal{L}_{q_j} ^{\rm T} \\
        \mathcal{L}_{v_j} ^{\rm T} 
    \end{bmatrix}
    = \begin{bmatrix}
        l_{q_j} ^{\rm T} \\
        l_{v_j} ^{\rm T} 
    \end{bmatrix}
    + \begin{bmatrix}
        \delta_{q_{j}} ^{\rm T} (q_{j}, q_{j+1}) & O \\
        O & I_n 
    \end{bmatrix}
    \begin{bmatrix}
        \lambda_{j+1} \\
        \gamma_{j+1}
    \end{bmatrix} \notag \\
    & + \begin{bmatrix}
        {\Gamma}_{q_j} ^{\rm T} ({y}_j) & {\bf v}_{q_j} ^{\rm T} ({z}_j) \\
        {\Gamma}_{v_j} ^{\rm T} ({y}_j) & {\bf v}_{v_j} ^{\rm T} ({z}_j)
    \end{bmatrix}
    \begin{bmatrix}
        \beta_j \\
        \mu_j 
    \end{bmatrix} 
    + \begin{bmatrix}
        \delta_{q_{j}} ^{\rm T} (q_{j-1}, q_j) \lambda_{j} \\
        - \gamma_{j}
    \end{bmatrix} = 0,
\end{align}
and
\begin{align}\label{eq:ldvLambdaj}
    \begin{bmatrix}
        \mathcal{L}_{\delta v_j} ^{\rm T} \\
        - \mathcal{L}_{\Lambda_j} ^{\rm T} 
    \end{bmatrix}
    = & \begin{bmatrix}
        l_{\delta v_j} ^{\rm T} \\
        - l_{\Lambda_j} ^{\rm T} 
    \end{bmatrix}
    + \begin{bmatrix}
        \gamma_{j+1} \\
        0 
    \end{bmatrix} 
    \notag \\ 
    & + \begin{bmatrix}
        {\Gamma}_{a_j} ^{\rm T} (\tilde{y}_j) & {\bf v}_{\delta v_j} ^{\rm T} (z_j) \\
        - {\Gamma}_{f_j} ^{\rm T} (\tilde{y}_j) & - {\bf v}_{f_j} ^{\rm T} (z_j)
    \end{bmatrix}
    \begin{bmatrix}
        \beta_j \\
        \mu_j
    \end{bmatrix} = 0,
\end{align}
where $\lambda_{j+1}, \allowbreak \gamma_{j+1}, \allowbreak \beta_j \in \mathbb{R}^n$, and $\mu_j \in \mathbb{R}^{n_f}$ are the Lagrange multipliers with respect to (\ref{eq:qvDiffImpulse}), (\ref{eq:inverseImpulseDynamics}), and (\ref{eq:v}), respectively.

These KKT conditions are also linearized into linear equations with respect to the Newton steps $\Delta q_i$, $\Delta v_i$, $\Delta \delta v_i$, $\Delta \Lambda_i$, $\Delta \lambda_{i+1}$, $\Delta \gamma_{i+1}$, $\Delta \beta_i$, and $\Delta \mu_i$, which are larger than those of the non-lifted counterpart with respect to $\Delta q_i$, $\Delta v_i$, $\Delta \lambda_{i+1}$, and $\Delta \gamma_{i+1}$.
Thus, we propose an efficient condensing method to reduce the size of the linear equation.
We first observe that the equality constraints (\ref{eq:inverseImpulseDynamics}) and (\ref{eq:v}) are linearized as 
\begin{align}\label{eq:impulseDynamicsLinearized}
    & \begin{bmatrix}
        M (q_j) & J ^{\rm T} (q_j) \\
        J (q_j) & O 
    \end{bmatrix}
    \begin{bmatrix}
        \Delta \delta v_j \\
        - \Delta \Lambda_j 
    \end{bmatrix} 
    \notag \\ 
    & = - \begin{bmatrix}
        {\Gamma}_{q_j} (y_j) & {\Gamma}_{v_j} (y_j) \\
        {\bf v}_{q_j} (z_j) & {\bf v}_{v_j} (z_j)
    \end{bmatrix}
    \begin{bmatrix}
        \Delta q_j \\
        \Delta v_j 
    \end{bmatrix}
    - \begin{bmatrix}
        {\Gamma} (y_j) \\
        {\bf v} (z_j)
    \end{bmatrix}.
\end{align}
Therefore, we can express $\Delta \delta v_j$ and $\Delta \Lambda_j$ with the linear combinations of $\Delta q_j$ and $\Delta v_j$ if we compute (\ref{equ:MJtJinv}) for $q_j$.
Next, we observe that (\ref{eq:ldvLambdaj}) is linearized into 
\begin{align}\label{eq:ldvLambdajLin}
    & \begin{bmatrix}
        \mathcal{L}_{\delta v_j} ^{\rm T} \\
        - \mathcal{L}_{\Lambda_j} ^{\rm T} 
    \end{bmatrix}
    + \begin{bmatrix}
        \mathcal{L}_{\delta v_j y_j} \\
        - \mathcal{L}_{\Lambda_j y_j} 
    \end{bmatrix}
    \Delta y_j
    + \begin{bmatrix}
        \Delta \gamma_{j+1} \\
        0 
    \end{bmatrix}
    \notag \\
    & + \begin{bmatrix}
        M (q_j) & J ^{\rm T} (q_j) \\
        J (q_j) & O
    \end{bmatrix}
    \begin{bmatrix}
        \Delta \beta_j \\
        \Delta \mu_j 
    \end{bmatrix} = 0.
\end{align}
Therefore, we can also express the Newton steps of the dual variables $\Delta \beta_j$ and $\Delta \mu_j$ with the linear combinations of $\Delta q_j$, $\Delta v_j$, and $\Delta \gamma_{j+1}$ if we solve (\ref{equ:MJtJinv}) for $q_j$.
Subsequently, the linear equations to be solved in the Newton-type iterations are reduced with respect to $\Delta q_j$, $\Delta v_j$, $\Delta \lambda_{j+1}$, and $\Delta \gamma_{j+1}$ and are expressed as 
\begin{align}
    & \begin{bmatrix}
        \tilde{F}_{q, j} \\ 
        \tilde{F}_{v, j}
    \end{bmatrix} 
    + \begin{bmatrix}
        \delta_{q_{j}} (q_{j}, q_{j+1}) & O \\ 
        \tilde{F}_{vq, j} & \tilde{F}_{vv, j}
    \end{bmatrix} 
    \begin{bmatrix}
        \Delta q_{j} \\ 
        \Delta v_{j}
    \end{bmatrix} 
    \notag \\
    &
    + \begin{bmatrix}
        \delta_{q_{j+1}} (q_{j}, q_{j+1}) & O \\ 
        O & - I
    \end{bmatrix}
    \begin{bmatrix}
        \Delta q_{j+1} \\ 
        \Delta v_{j+1}
    \end{bmatrix} = 0
\end{align}
and
\begin{align}\label{eq:lqviConImpulse}
    & \begin{bmatrix}
        \tilde{\mathcal{L}}_{q_j} ^{\rm T} \\
        \tilde{\mathcal{L}}_{v_j} ^{\rm T} 
    \end{bmatrix}
    + \begin{bmatrix}
        \tilde{\mathcal{L}}_{q_j q_j} &
        \tilde{\mathcal{L}}_{q_j v_j} \\
        \tilde{\mathcal{L}}_{v_j q_j} &
        \tilde{\mathcal{L}}_{v_j v_j} 
    \end{bmatrix}
    \begin{bmatrix}
        \Delta q_j \\
        \Delta v_i
    \end{bmatrix} 
    + \begin{bmatrix}
        \delta_{q_{j}} ^{\rm T} (q_{j-1}, q_j) \Delta \lambda_{j} \\
        - \Delta \gamma_{j}
    \end{bmatrix} 
    \notag \\
    & + \begin{bmatrix}
        \delta_{q_{j}} ^{\rm T} (q_{j}, q_{j+1}) 
        & \tilde{F}_{vq, j} \\
        O 
        & \tilde{F}_{vv, j}
    \end{bmatrix}
    \begin{bmatrix}
        \Delta \lambda_{j+1} \\
        \Delta \gamma_{j+1}
    \end{bmatrix} 
    = 0,
\end{align}
for which we omit the definitions of the terms with a superscript tilde because they are similar to those of the lifted contact dynamics.
After solving the linear equations and obtaining $\Delta q_j$, $\Delta v_j$, $\Delta \lambda_{j+1}$, and $\Delta \gamma_{j+1}$, we can easily compute $\Delta \delta v_j$, $\Delta \Lambda_j$, $\Delta \beta_j$, and $\Delta \mu_j$ using (\ref{eq:impulseDynamicsLinearized}) and (\ref{eq:ldvLambdajLin}), which is the expansion procedure of the impulse stage.

\subsection{Riccati recursion}
After performing the above condensing, the Newton step computation is reduced to solving a set of linear equations with respect to $\Delta q_i$, $\Delta v_i$, $\Delta u_i$, $\Delta \lambda_i$, and $\Delta \gamma_i$.
This can be seen as the KKT conditions of an LQR subproblem with the state $[\Delta q_i ^{\rm T} \;\; \Delta v_i ^{\rm T}]^{\rm T}$ and control input $\Delta u_i$. 
That is, solving the linear equations after condensing is equivalent to solving the LQR subproblem. 
We then utilize the Riccati recursion algorithm \cite{bib:mpcbook}, which can solve the LQR subproblem only with $O(N)$ complexity while the direct inversion of the KKT matrix requires an $O(N^3)$ computational burden. 
For example, DDP uses the Riccati recursion for single-shooting OCPs with a nonlinear forward pass and can efficiently solve very large-scale problems \cite{bib:legged2, bib:DDPContact, bib:crocoddyl}.

\subsection{Primal-dual interior-point method}
OCPs of rigid-body systems can contain many inequality constraints, e.g., joint angle limits, joint angular velocity limits, joint torque limits, and friction cones.
To treat such a large number of inequality constraints including nonlinear ones with the Riccati recursion algorithm efficiently, we use the primal-dual interior-point method \cite{bib:nocedal} (PDIPM).
In the PDIPM, only certain terms related to the second- and first-order derivatives of the logarithmic barrier functions are added to the Hessians $\mathcal{L}_{y_i y_i}$ and residuals $\mathcal{L}_{y_i}$, respectively, and there are no effects on the proposed condensing method and the LQR subproblem. 
We can then apply the Riccati recursion algorithm regardless of the inequality constraints. 
After solving the LQR subproblem and obtaining $\Delta y_i$, we can efficiently compute the Newton steps of the slack variables and Lagrange multipliers of the inequality constraints.

\subsection{Algorithm}
We summarize the single Newton iteration, that is, the computation of the Newton steps of the proposed method, in Algorithm 1.
First, we compute the Hessians of the Lagrangian, Jacobians of the constraints, and residuals in the KKT conditions (line 2).
We utilize the Gauss-Newton method \cite{bib:mpcbook} to avoid computing the second-order derivatives of the dynamics.
The modification of the Hessians and KKT residuals due to the PDIPM is also performed in this step.
Second, we compute the matrix inversion (\ref{equ:MJtJinv}) and form the condensed Hessians, Jacobians, and KKT residuals (lines 3 and 4).
These two steps are fully parallelizable because of the multiple-shooting formulation.
We then solve the LQR subproblem, e.g., using the Riccati recursion (line 6).
Finally, we compute the condensed Newton steps from the solution of the LQR subproblem (line 8).
The Newton steps of the slack variables and Lagrange multipliers of the PDIPM are also computed in this step.

\begin{algorithm}[tb]
\caption{Newton step computation in the lifted contact dynamics scheme}
\label{alg1}
\begin{algorithmic}[1]
    \Require Initial state ${x} (t_0)$, the current iterate $y_0, ..., y_{N-1}$, $q_N$, $v_N$, $\left\{ u_i \right\}_{i \in \mathcal{I}}$, $\lambda_0, ..., \lambda_N$, $\gamma_0, ..., \gamma_N$, $\beta_0, ..., \beta_{N-1}$, $\mu_0, ..., \mu_{N-1}$, and $\left\{ \nu_i \right\}_{i \in \mathcal{I}}$
    \Ensure Newton steps $\Delta y_0, ..., \Delta y_{N-1}$, $\Delta q_{N}$, $\Delta v_N$, $\left\{ \Delta u_i \right\}_{i \in \mathcal{I}}$, $\Delta \lambda_0, \allowbreak ..., \allowbreak \Delta \lambda_N$, $\Delta \gamma_0, \allowbreak ..., \allowbreak \Delta \gamma_N$, $\Delta \beta_0, \allowbreak ..., \allowbreak \Delta \beta_{N-1}$, $\Delta \mu_0, \allowbreak ..., \allowbreak \Delta \mu_{N-1}$, and $\left\{ \Delta \nu_i \right\}_{i \in \mathcal{I}}$
    \For{$i=0,\cdots,N$} {\bf in parallel} 
    \State Compute the Hessian of the Lagrangian (e.g., $\mathcal{L}_{q_i q_i}$), Jacobians of the constraints (e.g., ${\rm ID}_{q_i}$ and ${\bf a}_{q_i}$), and residuals in the KKT conditions (e.g., (\ref{eq:qvDiff}) and $\mathcal{L}_{q_i}$).
    \State Compute the matrix inversion (\ref{equ:MJtJinv}).
    \State Form the condensed Hessians, Jacobians, and KKT residuals (e.g., $\tilde{\mathcal{L}}_{q_i q_i}$, $\tilde{F}_{qq, i}$, and $\tilde{\mathcal{L}}_{q_i}$).
    \EndFor
    \State Compute $\Delta q_0, ..., \Delta q_N$, $\Delta v_0, ..., \Delta v_N$, $\left\{ \Delta u_i \right\}_{i \in \mathcal{I}}$, $\Delta \lambda_0, ..., \Delta \lambda_N$, and $\Delta \gamma_0, ..., \Delta \gamma_N$ by solving the LQR subproblem, e.g., using the Riccati recursion.
    \For{$i=0,\cdots,N-1$} {\bf in parallel} 
    \State Compute the condensed Newton steps ($\Delta a_i$, $\Delta f_i$, $\Delta \beta_i$, $\Delta \mu_i$, and $\Delta \nu_i$ for $i \in \mathcal{I}$ or $\Delta \delta v_j$, $\Delta \Lambda_j$, $\Delta \beta_j$, and $\Delta \mu_j$ for $j \in \mathcal{J}$).
    \EndFor
\end{algorithmic}
\end{algorithm}

\subsection{Comparison with existing methods}
Here, we re-summarize the comparison between the proposed method and the existing methods:
\subsubsection{Comparison with non-lifted formulations \cite{bib:humanLikeRunnig, bib:DDPContact, bib:crocoddyl, bib:DDP:jumpRobot}}
The non-lifted formulations \cite{bib:humanLikeRunnig, bib:DDPContact, bib:crocoddyl, bib:DDP:jumpRobot} eliminate $a_i$, $f_i$, $\delta v_i$, and $\Lambda_i$ as (nonlinear) functions of $q_i$, $v_i$, and $u_i$.
Therefore, the costs and constraints including $a_i$, $f_i$, $\delta v_i$, and $\Lambda_i$ (e.g., the friction cone constraints) can be highly nonlinear.
The proposed method regard $a_i$, $f_i$, $\delta v_i$ as decision variables to alleviate such high nonlinearities.
The computational costs per Newton-type iteration of the proposed method and these methods are similar owing to the proposed efficient condensing procedure.
\subsubsection{Comparison with inverse dynamics-based algorithm \cite{bib:InvDynOCP}}
The inverse dynamics-based method of \cite{bib:InvDynOCP} regards $a_i$ and $f_i$ as the optimization variables: the convergence behavior is identical to that of the proposed method.
However, the proposed condensing method can be more efficient than that presented in \cite{bib:InvDynOCP}.
The proposed method requires $O (n_a^3)$ computational burden to compute the Newton steps with the Riccati recursion.
The method of \cite{bib:InvDynOCP} requires $O ((n+n_f)^3)$, which is inefficient when the numbers of contacts and passive joints (i.e., $n_f$ and $n - n_a$) are large.

\section{Numerical Experiments: Whole-Body Optimal Control of Quadrupedal Gaits}\label{section:expriment}
\subsection{Experimental settings}

\begin{table}[tb]
\centering
\caption{Settings of the OCPs for the five quadrupedal gaits}
\begin{tabular}{l|lllll}
          & Walking & Trotting & Pacing & Bounding & Jumping \\ \hline 
$N$       & 107     & 57       & 71     & 71       & 71 \\
$\Delta \tau$ & 0.02    & 0.02     & 0.01   & 0.01 & 0.01 
\end{tabular}
\end{table}

To demonstrate the effectiveness of the proposed lifting approach over the conventional non-lifting approaches, we conducted numerical experiments on the whole-body optimal control of quadrupedal robot ANYmal for five gaits, that is, walking, trotting, pacing, bounding, and jumping, subject to the friction cone constraints.
We considered the polyhedral-approximated friction cone constraint for each contact force expressed in the world frame 
$[f_x \;\, f_y \;\, f_z]$ as
\begin{equation}\label{eq:frictionCone}
    \begin{bmatrix}
        f_x + \frac{\mu}{\sqrt{2}} f_z  \\
        - f_x + \frac{\mu}{\sqrt{2}} f_z  \\
        f_y + \frac{\mu}{\sqrt{2}} f_z  \\
        - f_y + \frac{\mu}{\sqrt{2}} f_z  \\
        f_z
    \end{bmatrix} \geq 0, 
\end{equation}
where $\mu > 0$ is the friction coefficient, and we set it as $0.7$.
The problem settings of the OCPs for the five gaits are shown in Table I.
We compared the following five methods:
\begin{itemize}
    \item DMS-LCD: DMS based on the proposed lifted contact dynamics 
    \item DMS-CD: DMS based on the conventional non-lifted contact dynamics (e.g., used in \cite{bib:humanLikeRunnig})
    \item DMS-ID: Inverse dynamics-based DMS algorithm \cite{bib:InvDynOCP} 
    \item FDDP: Feasibility-driven DDP (FDDP), a variant of DDP with Gauss-Newton method that improves the numerical robustness \cite{bib:crocoddyl}, based on the conventional non-lifted contact dynamics 
    \item iLQR \cite{bib:ilqr}: DDP with Gauss-Newton method based on the conventional non-lifted contact dynamics
\end{itemize}
The following  were compared in this study: 
1) The convergence property between the lifted and non-lifted formulations (DMS-LCD vs. DMS-CD);
2) The CPU time per Newton iteration between the two condensing methods for the lifted formulation (DMS-LCD vs. DMS-ID); 
3) The proposed method with a state-of-the-art OCP solver \cite{bib:crocoddyl} (DMS-LCD vs. FDDP and iLQR).
DMS-LCD, DMS-CD, and DMS-ID used Gauss--Newton Hessian approximation, the PDIPM method for inequality constraints, and the Riccati recursion to solve the LQR subproblems.
FDDP and iLQR used the relaxed barrier function method (ReB) \cite{bib:relaxedBarrier}, which is a popular constraint-handling approach in DDP-type methods (e.g., used in \cite{bib:DDP:jumpRobot}) that enables the iteration to lie outside the feasible region.
Note that the ReB and the PDIPM correspond to the same barrier function under the same barrier parameter and feasible solution \cite{bib:nocedal, bib:relaxedBarrier}.
That is, the five methods consider the same barrier functions for the inequality constraints.
We implemented DMS-LCD\footnote{Our open-source implementation of DMS-LCD is available online at \url{https://github.com/mayataka/robotoc}}, DMS-ID, and DMS-CD in C++ and used \texttt{Pinocchio} \cite{bib:pinocchio}, an efficient C++ library used for rigid-body dynamics and its analytical derivatives, to compute the dynamics and its derivatives of the quadrupedal robot.
We used OpenMP for parallel computing (e.g., lines 1--5 of Algorithm 1). 
For FDDP and iLQR, we used the \texttt{Crocoddyl} \cite{bib:crocoddyl} framework, which was also implemented in C++, \texttt{Pinocchio} for rigid-body dynamics, and OpenMP for parallel computing.
FDDP and iLQR consider the Baumgarte's stabilization method of the form in (\ref{eq:contactDynamics}) as well as DMS-based methods. 

We designed $V_f (\cdot)$ and $l_i (\cdot)$ in (\ref{eq:cost}) as least-square objectives to track the pre-defined feet and center of mass trajectories.
We set the stabilization parameters in (\ref{eq:a}) as $\alpha = \beta = 25$.
We approximated the contact position constraints (\ref{eq:p}) using a quadratic penalty function for simplicity.
DMS-LCD, DMS-CD, and DMS-ID did not use any regularization on the Hessian and used only the fraction-to-boundary-rule \cite{bib:nocedal} for the step-size selection.
FDDP and iLQR used an adaptive regularization on the Hessian and a line-search based on the Goldstein condition in the step-size selection, as provided in the \texttt{Crocoddyl} solver. 
We conducted all the experiments for the two barrier parameters $\epsilon$: $1.0 \times 10^{-1}$ and $1.0 \times 10^{-4}$, which represent large and small barrier parameters, respectively.
$\epsilon$ is fixed during each experiment, which corresponds to suboptimal MPC cases.
We set the relaxation parameter of the ReB for each $\epsilon$ and each gait as large as possible while keeping the optimal solution satisfying (\ref{eq:frictionCone}), which means that the ``ideal'' ReB (that is, a favorable experimental setting for ReB-based FDDP and iLQR) was considered in the experiments 
The experiment was conducted with an octa-core CPU Intel Core i9-9900 @3.10 GHz, and all the algorithms were compiled using the GCC compiler with \texttt{-O3 -DNDEBUG -march=native} options.
We used eight threads in parallel computing.

\subsection{Results and discussion}
\begin{figure}[tb]
\centering
\includegraphics[scale=0.57]{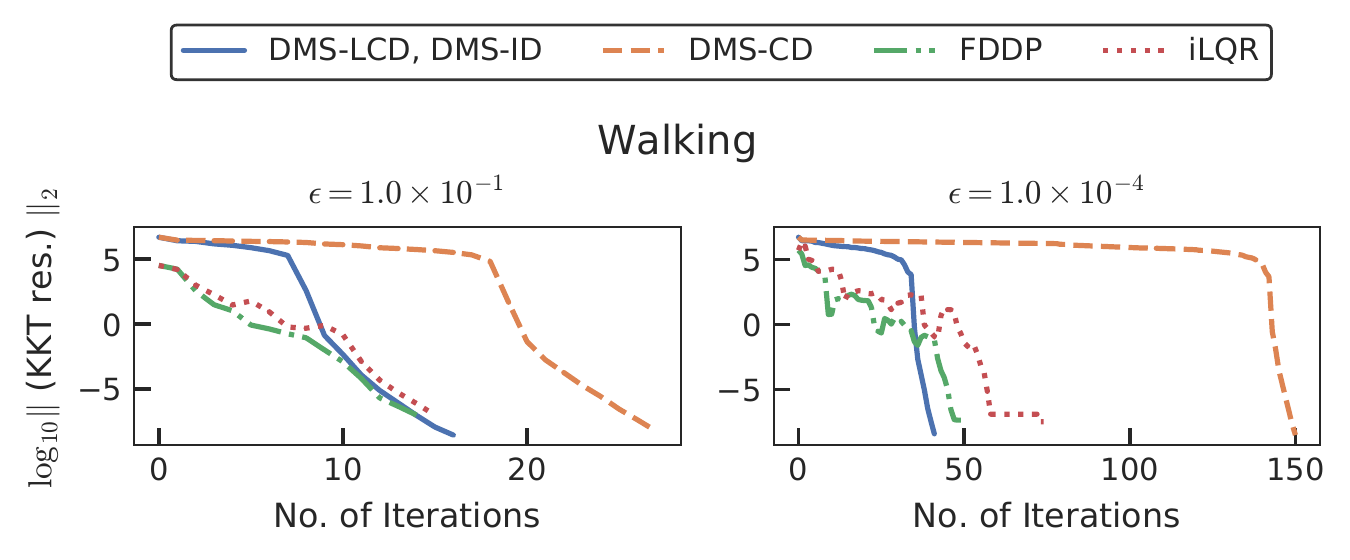}
\includegraphics[scale=0.57]{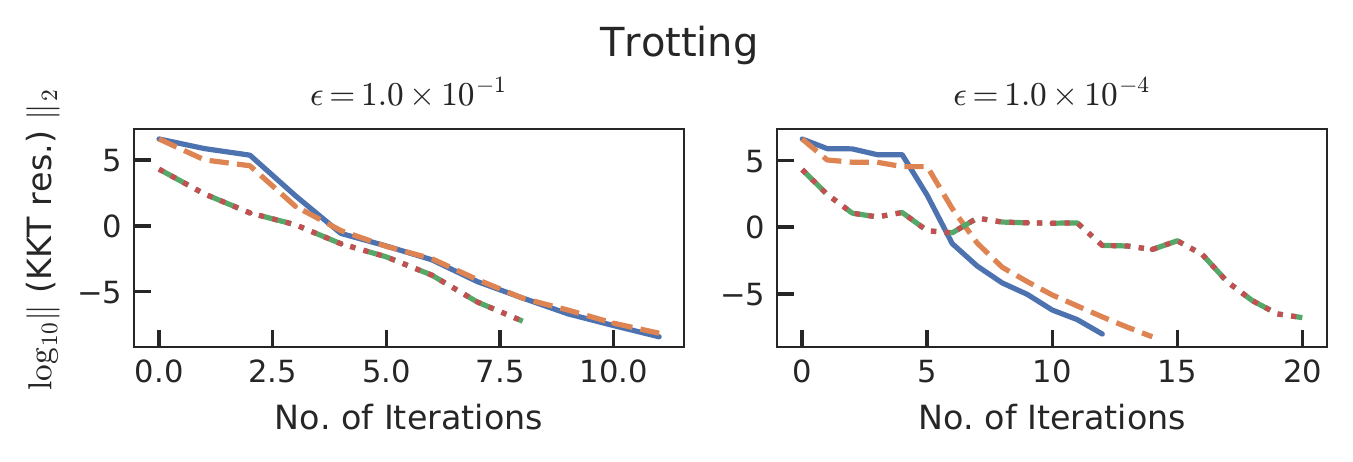}
\includegraphics[scale=0.57]{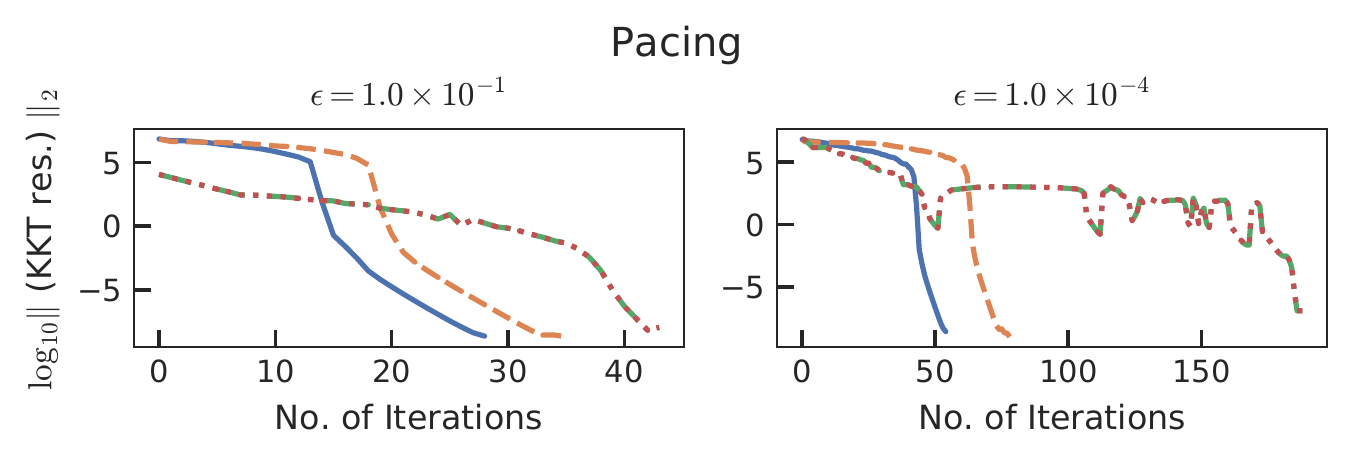}
\includegraphics[scale=0.57]{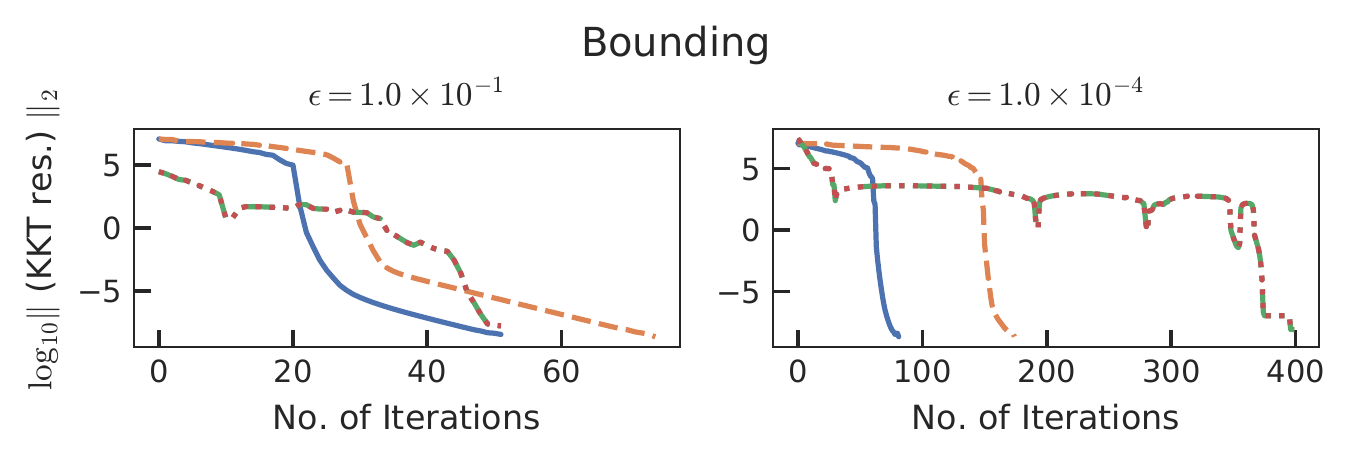}
\includegraphics[scale=0.57]{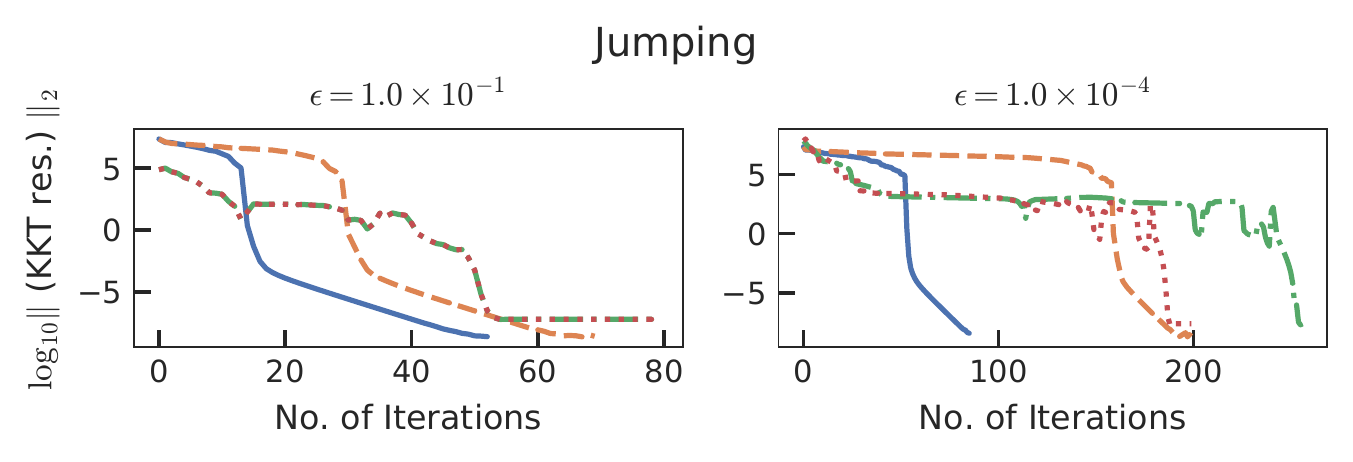}
\caption{$\log_{10}$ scaled $l_2$-norms of the KKT residuals (KKT res. in short) of DMS-LCD and DMS-ID (solid lines), DMS-CD (dashed lines), FDDP (dash-dotted lines), and iLQR (dotted lines) over the five quadrupedal gaits subject to the friction cone constraints considered in the interior-point methods with the fixed barrier parameters $\epsilon = 1.0 \times 10 ^{-1}$ and $\epsilon = 1.0 \times 10 ^{-4}$.}
\label{fig:KKT}
\end{figure}

Fig. \ref{fig:KKT} shows $\log_{10}$ scaled $l_2$-norms of the KKT residuals with respect to the number of iterations of these five methods over the five quadrupedal gaits with the barrier parameters $\epsilon = 1.0 \times 10 ^{-1}$ and $\epsilon = 1.0 \times 10 ^{-4}$.
Note that the convergence results of DMS-LCD and DMS-ID are illustrated in the same solid lines because they are identical.
As shown in Fig. \ref{fig:KKT}, the five methods resulted in a similar convergence when the barrier parameter was large ($\epsilon = 1.0 \times 10 ^{-1}$). 
In contrast, when the barrier parameter was small ($\epsilon = 1.0 \times 10 ^{-4}$), the proposed DMS-LCD converged significantly faster than the other methods, particularly in the aggressive gaits (pacing, bounding, and jumping).
This is because, in the interior-point methods, a smaller barrier parameter (i.e., a more accurate solution) corresponds to a slower and more difficult convergence due to the high nonlinearity \cite{bib:nocedal}.
Nevertheless, the lifting methods (DMS-LCD and DMS-ID) achieved fast convergence even with the small barrier parameter.
The faster convergence of DMS-LCD than DMS-CD particularly shows the effectiveness of the lifted formulation.

\begin{figure}[tb]
\centering
\includegraphics[scale=0.57]{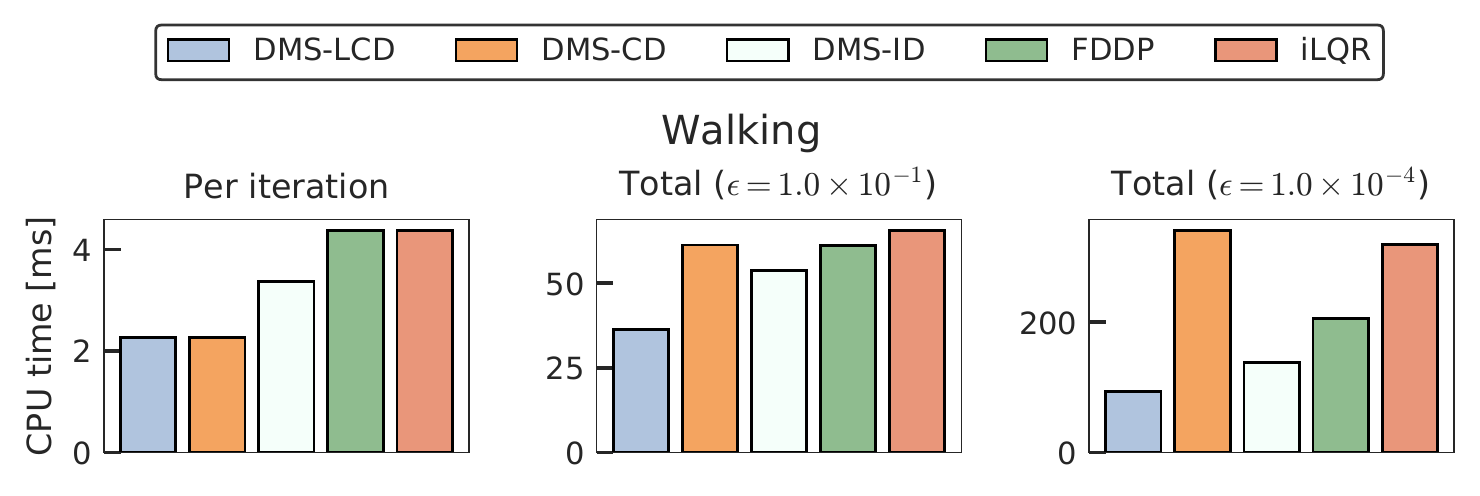}
\includegraphics[scale=0.57]{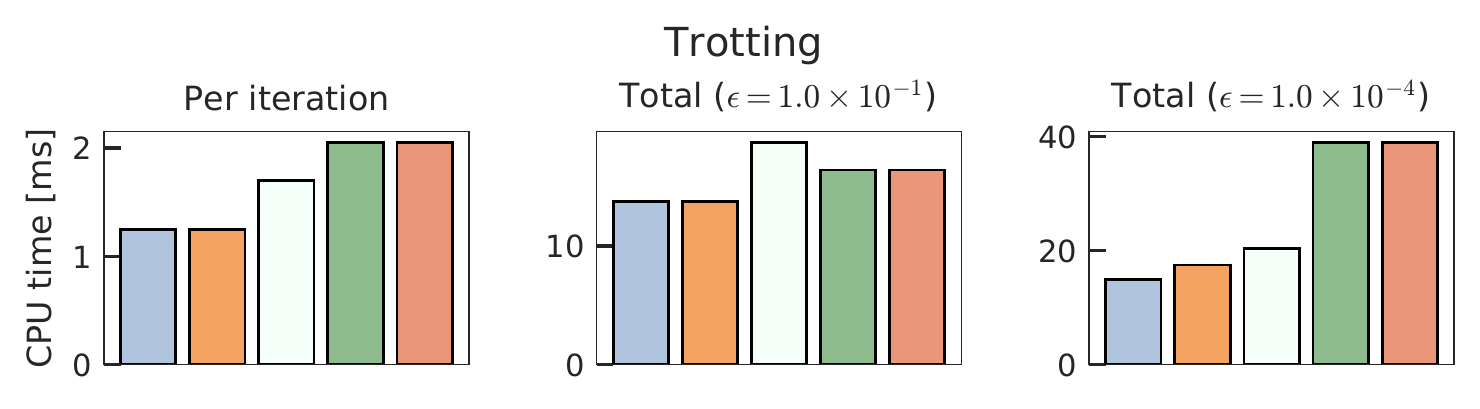}
\includegraphics[scale=0.57]{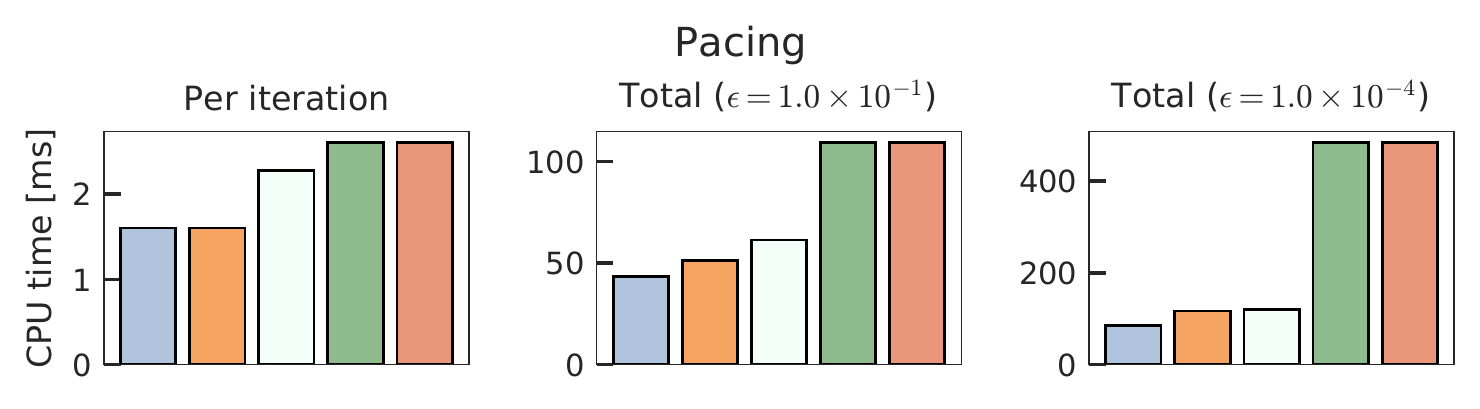}
\includegraphics[scale=0.57]{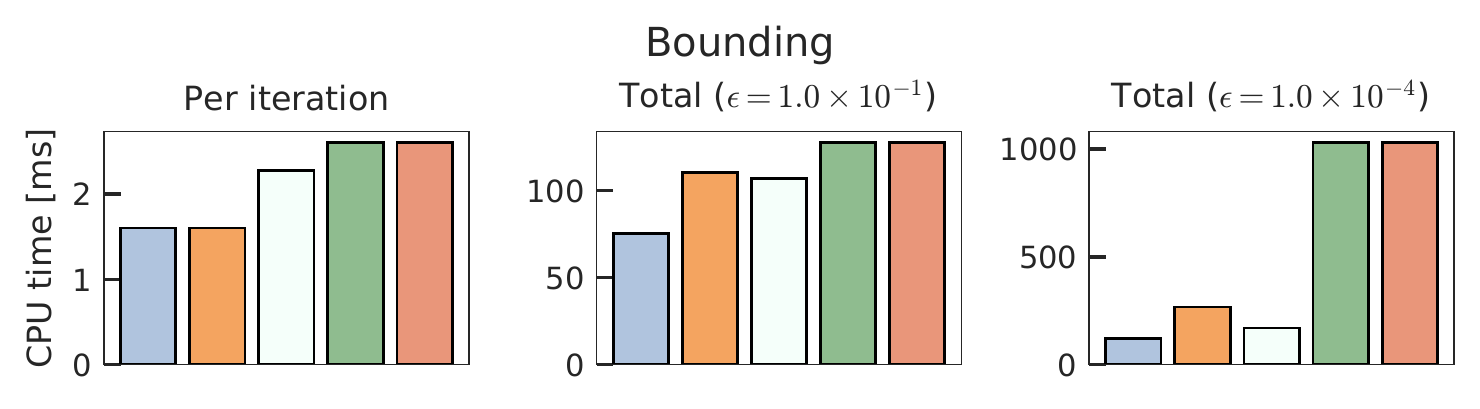}
\includegraphics[scale=0.57]{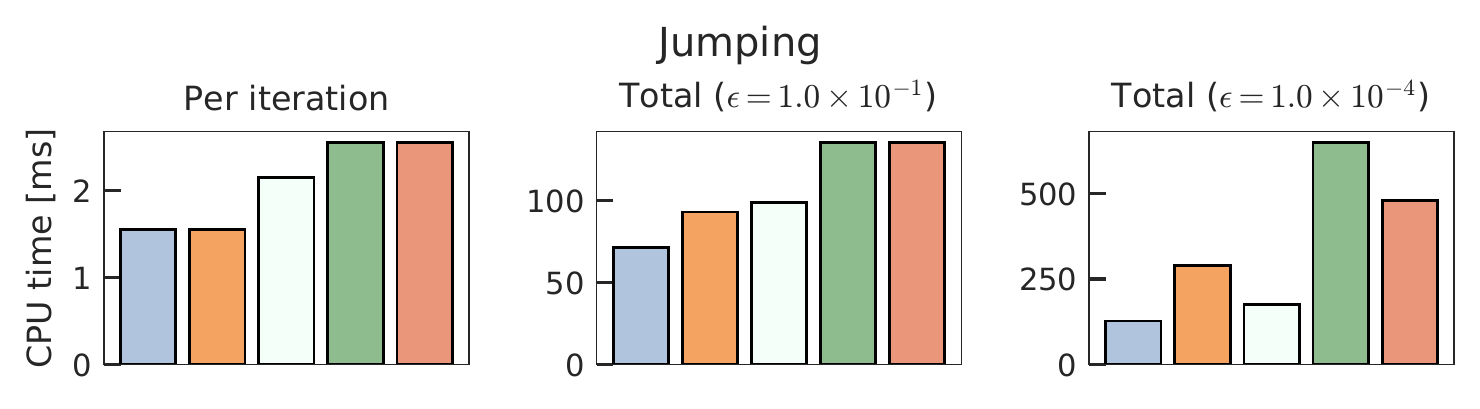}
\caption{CPU time per iteration and the total CPU time until convergence of DMS-LCD, DMS-CD, DMS-ID, FDDP, and iLQR over the five quadrupedal gaits subject to the friction cone constraints considered in the interior-point methods with the fixed barrier parameters $\epsilon = 1.0 \times 10 ^{-1}$ and $\epsilon = 1.0 \times 10 ^{-4}$.
}
\label{fig:CPUTime}
\end{figure}

Fig. \ref{fig:CPUTime} shows the CPU time per iteration and the total CPU time until convergence of the five methods.
First, we observed that the CPU times per iteration of DMS-LCD and DMS-CD were almost identical and around 1.5 times faster than those of DMS-ID, which shows the efficiency of the proposed condensing algorithm over the previous one \cite{bib:InvDynOCP}. 
Furthermore, DMS-LCD and DMS-CD were 1.6 to 2 times faster than those of FDDP and iLQR because DMS could leverage parallel computing in the computation of the KKT residual, whereas single shooting methods such as FDDP and iLQR were required to compute the contact and impulse dynamics (\ref{eq:contactDynamics}) and (\ref{eq:impulseDynamics}) over the horizon with a single thread.
We compared the total computational time, and DMS-LCD had the fastest computational time in all the scenarios.
It was more than twice as fast as the other non-lifted methods in several scenarios.

\begin{figure}[tb]
\centering
\includegraphics[scale=0.65]{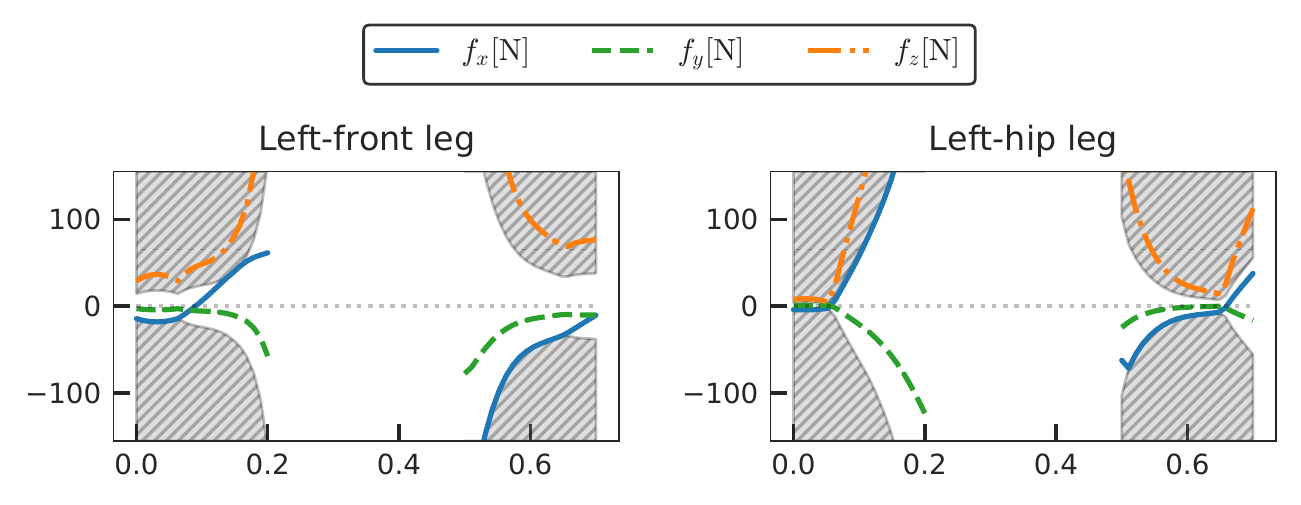}
\caption{Time histories of the contact forces expressed in the world frame $[f_x \;\, f_y \;\, f_z]$ of left legs in the jumping motion. The infeasible regions of $f_x$ and $f_y$ due to the friction cone constraints (\ref{eq:frictionCone}) are the filled gray hatches. 
The infeasible region of $f_z \geq 0$ is the lower-half space of each plot.}
\label{fig:force}
\end{figure}

Fig. \ref{fig:force} shows the time histories of the contact forces of the left feet in the motion.
It shows that the friction cone constraints (\ref{eq:frictionCone}) were active and satisfied during the jumping.
A supplemental video including all the five gaits is available at \url{https://youtu.be/jb7gGnblQ7s}. 

We further conducted several experiments with the various Baumgarte's weight parameters $\alpha = \beta > 0$ in (\ref{eq:a}), the details of which are omitted here owing to the space limitation.
We observed that as the weight parameter increased, the proposed method had greater advantages in terms of the convergence speed.
In contrast, if the weight parameter was small, the proposed method had no such advantages. 
Note that a larger weight parameter typically results in a more accurate solution in terms of the original contact position constraints (\ref{eq:p}).
This observation confirms that the proposed method can address high nonlinearity because the large Baumgarte's weight parameters cause high nonlinearity \cite{bib:baumgarteParameters}.

\section{Conclusions}\label{section:conclu}

We proposed a novel lifting approach for optimal control of rigid-body systems with contacts to improve the convergence properties of Newton-type methods.
To relax the high nonlinearity, we considered the acceleration and contact forces as the optimization variables and the inverse dynamics and acceleration constraints of the contact frames as equality constraints.
We eliminated the update of these additional variables from the linear equation for Newton-type method in an efficient manner. 
As a result, the computational cost per Newton-type iteration is almost identical to that of the conventional non-lifted one that embeds contact dynamics in the state equation. 
We conducted numerical experiments on the whole-body optimal control of various quadrupedal gaits subject to the friction cone constraints considered in the interior-point methods and demonstrated that the proposed method can significantly increase the convergence speed to more than twice that of the conventional non-lifted approaches.
A future research direction is to combine the proposed method with the switching-time optimization methods \cite{bib:DDP:jumpRobot} to study MPC for legged robots.








\section*{Acknowledgment}
This work was partly supported by JST SPRING, Grant Number JPMJSP2110, and JSPS KAKENHI, Grant Numbers JP22J11441 and JP22H01510.


\bibliographystyle{IEEEtran}
\bibliography{IEEEabrv, ieee}


\end{document}